\documentclass[11pt,a4paper]{article}
\usepackage{beamerarticle}

\mode<article>{\usepackage{fullpage}}
\mode<presentation>{\usetheme{Berlin}}

\usepackage{amsfonts, amsmath}
\usepackage{graphicx}

\usepackage{color}

\usepackage[english]{babel}
\usepackage{pgf}
\usepackage[utf8]{inputenc}
\usepackage{fontenc}

\newcommand{\mc}[1]{\mathcal{#1}}

\newcommand{\com}[1]{\opt{draft}{\textcolor{red}{
$\LHD$ #1 $\RHD$\marginpar{\textcolor{red}{$\begin{lema}acksquare$}}}}}

\newcommand{\comb}[1]{\opt{draft}{\textcolor{blue}{
$\LHD$ #1 $\RHD$\marginpar{\textcolor{blue}{$\begin{lema}acksquare$}}}}}

\def\qed{\hfill {\large ${\sqcup\!\!\!\!\sqcap}$}}

\newenvironment{demo}{{\bf Proof }}
{\qed \\}

\newcommand{\grad}{{\rm grad}}
\newcommand{\dv}{{\rm div}}

\newcommand{\re}{\mathbb R}

\newcommand{\<}{\left<}
\renewcommand{\(}{\left(}
\newcommand{\lb}{\label}
\newcommand{\nn}{\nonumber}

\newcommand{\fracc}{\displaystyle\frac}
\newcommand{\ds}{\displaystyle}
\renewcommand{\>}{\right>}
\renewcommand{\)}{\right)}
\newcommand{\flecha}{\longrightarrow}

\def\p{\varphi}

\def\vec{\overrightarrow}

\newcommand{\bde}{\begin{defi}}
\newcommand{\ede}{\end{defi}}

\newcommand{\be}{\begin{enumerate}}
\newcommand{\ee}{\end{enumerate}}

\newcommand{\ba}{\begin{array}}
\newcommand{\ea}{\end{array}}

\def\og{{\overline g}}

\def\oM{{\overline M}}
\def\ona{{\overline \nabla}}

\def\oDelta{\overline{\Delta}}

\def\odv{{\overline{\rm div}}}

\def\wal{{\widetilde \alpha }}

\def\wps{{\widetilde \psi}}
\def\wHp{{\widetilde{\vec{H}}_{\!\!\!{\wps}}}}

\def\wvH{{\widetilde{\vec{H}}}}

\def\vHp{{\vec{H}}_{\!\!\!{\psi}}}

\newtheorem{defi}{\hspace{12pt} Definition}
\newtheorem{teor}{\hspace{12pt} Theorem}
\newtheorem{prop}[teor]{\hspace{12pt} Proposition}

\newtheorem{lema}[teor]{\hspace{12pt} Lemma}
\newtheorem*{lema*}{\hspace{12pt} Lemma}

\newtheorem{nota}{\hspace{12pt} Remark}

\newcommand{\ben}{\begin{enumerate}}
\newcommand{\een}{\end{enumerate}}
\newcommand{\bi}{\begin{itemize}}
\newcommand{\ei}{\end{itemize}}
\newcommand{\bec}{\begin{equation}}
\newcommand{\eec}{\end{equation}}
\newcommand{\beca}{\begin{equation*}}
\newcommand{\eeca}{\end{equation*}}
\newcommand{\bal}{\begin{align}}
\newcommand{\aal}{\end{align}}
\newcommand{\bala}{\begin{align*}}
\newcommand{\aala}{\end{align*}}

\newcommand{\Ei}{\operatorname{Ei}}
\newcommand{\erfi}{\operatorname{erfi}}

\begin{document}

\title{Reilly's type inequality for the Laplacian associated to a density related with shrinkers for MCF}

\author{ M. Carmen Domingo-Juan, Vicente Miquel  and Jonathan J. Zhu
\footnote{Mathematics Subject
Classification(2010) 53C42, 52C21}}

\date{}

\maketitle

\vspace{-1cm}
\begin{abstract}  
Let $(\oM,\<,\>,e^\psi)$ be a Riemannian manifold with density, and let $M$ be a closed $n$-dimensional submanifold  of $\oM$ with the induced metric and density. We give an upper bound on the first eigenvalue $\lambda_1$ of the closed eigenvalue problem for $\Delta_\psi$ (the Laplacian on $M$ associated to the density) in terms of the average of the norm of the vector $\vHp+ \ona\psi$ with respect to the volume form induced by the density, where $\vHp$ is the mean curvature of $M$ associated to the density $e^\psi$. \\
When $\oM=\re^{n+k}$ or $\oM=S^{n+k-1}$, the equality between $\lambda_1$ and its bound implies that $e^\psi$ is a Gaussian density ($\psi(x) = \frac{C}{2} |x|^2$, $C<0$), and $M$ is a shrinker for the mean curvature flow (MCF) on $\re^{n+k}$. \\
We prove also that  $\lambda_1 =-C$ on the standard shrinker torus of revolution. \\
Based on this and on the Yau's conjecture on the first eigenvalue of minimal submanifolds of $S^n$, we conjecture that the equality $\lambda_1=-C$ is true for all the shrinkers of MCF in $R^{n+k}$.
\end{abstract}

\section{Introduction }\lb{In}

After the seminal works of Bleecker, Weiner (\cite{BW}) and  Reilly (\cite{Re}), the obtention of extrinsic upper bounds for the eigenvalues of the Laplacian $\Delta$  on a submanifold in the Euclidean Space or other more general ambient spaces has become an interesting and fruitful problem. As examples of these works, one can look  at the references in \cite{GMO} and \cite{CESG}. 

Especially relevant for us are the quoted work of Reilly and the paper \cite{he} by Heintze, where they got sharp bounds for the first eigenvalue of the Laplacian on a closed manifold embedded in Euclidean space, in the sphere or in spaces with sectional curvature bounded from above. They also prove that equality is attained only at minimal submanifolds of some sphere, in particular, at the sphere when the codimension of the submanifold is $1$ and the ambient space is Euclidean space.

In the last years  the study of the spectrum of the Laplacian associated to a density $\Delta_\psi$ has received an increasing interest. A sample is the works of Ma, Liu, Du (\cite{ma14,maliu,madu}), Cheng, Mejia, Zhu (\cite{chmezh12b}), Pigola, Rimoldi (\cite{piri}), Ding, Xin 
(\cite{dixi}), Du, Mao, Wang and Wu (\cite{DMWW}). Recently, Batista, Cavalcante and Pyo (\cite{BCP}) have studied Reilly's type inequalities for $\Delta_\psi$, generalizing to this operator the Heintze's results on submanifolds of Riemannian manifolds with curvature bounded from above.
 The Laplacian associated to a density has also been used in the study of Kähler-Einstein manifolds by Li (\cite{lilo}) and it appears in a natural way as a part of the operator used in \cite{comi} by Colding and Minicozzi in the study of the $F$-stability of shrinkers in mean curvature flow.

Our contribution to this line of research is: i)  a deeper understanding of the Reilly's type inequalities for the Laplacian associated to a density obtained in \cite{BCP}, with the observation that, in Euclidean space, equality holds only on shrinkers of the mean curvature flow (theorems \ref{TF} and \ref{CorF}), ii) the determination of the first eigenvalue of $\Delta_\psi$ on some Angenent's type torus, proving that it attains the bound given by the previous inequality (Theorem \ref{l1A}), and iii) the observation of a relation  between the problem of the characterization of the submanifolds that attain the bound with a conjecture of Yau on the first eigenvalue of the Laplacian in a closed submanifold of the sphere. 

Before stating our results, we introduce the basic concepts.

An $m$-dimensional Riemannian manifold $(\overline{M} , \overline{g})$ with density is a Riemannian manifold where volumes are measured with a weighted (smooth) function 
$
 e^{\psi} : \overline{M} \flecha \Bbb R 
$
in the following way. If $\Omega$ is a domain in $\overline{M}$, $M$ is a closed submanifold of dimension $n$, and  $\mu_\og$ and $\mu_g$ are their respective Riemannian volume elements, the volume elements associated to the density $e^\psi$ are $\overline\mu_\psi = e^{\psi} \ \mu_{\overline{g}}$ and $\mu_\psi = e^{\psi} \ \mu_{g}$ respectively. This gives, for the corresponding $m$-volume $V_{\psi} (\Omega)$ and $n$-volume $V_{\psi} (M)$,
\begin{equation}
V_{\psi} (\Omega) = \int_{\Omega} e^\psi \ \mu_{\overline{g}}, \qquad \quad 
V_{\psi} (M) = \int_{M} e^\psi \ \mu_{g} \label{VA}.
\end{equation}

The generalization of mean curvature vector for  manifolds with densities is denoted by $\vec{H_{\psi}}$ and given by 
\begin{equation}
\vec{H_{\psi}} = \vec H -  \(\overline{\nabla} \psi\)^\bot , \label{curvden}
\end{equation}
where $\(\overline{\nabla} \psi\)^\bot$ denotes the component of $\overline{\nabla} \psi$ orthogonal to $M$, and we define the mean curvature $\vec{H}$ as the trace of the second fundamental form. With this definition $\vec{H_{\psi}}$ is the $L^2$-gradient (with respect to the measure $\mu_\psi$) of the functional $V_\psi$ defined on $n$-dimensional submanifolds of $\oM$. 

When $\psi(x)=a + \frac12\ C\ |x|^2$ (the Gaussian density), it is well known (for instance see \cite{sm12}) that the submanifolds $F:M\flecha \re^m$ satisfying $\vec{H_{\psi}}=0$, that is
$\overrightarrow{H} = C \ F^{\perp}$ (where $^\perp$ denotes the component orthogonal to $M$), are the self-similar solutions of the mean curvature flow. They can be: {\it shrinkers}, when $C<0$, which contract to a point under the flow, {\it minimal submanifolds}, when $C=0$, and {\it expanders}, when $C>0$, which enlarge under the flow.

Along this paper we shall use  $\odv$,  $\oDelta$ and $\ona$ to denote the divergence, Laplacian and gradient or covariant derivative respectively in the ambient manifold $\oM$, and $\dv$,  $\Delta$, $\nabla$ to denote the corresponding operators on the submanifold $M$.

In a manifold with density, the divergence and Laplacian associated to the volume form $\mu_{\psi} $  are defined as
\begin{equation}
\mathrm{div}_{\psi} X \mu_{\psi} = \mathcal{L}_X \mu_{\psi}, \qquad \quad 
\Delta _{\psi} f = \mathrm{div} _{\psi} \mathrm{grad} \ f ,
\label{divlap}
\end{equation}
and called $\psi$-divergence and $\psi$-Laplacian (drifted or drifting Laplacian in some references). Our first result is  a remark on the   Reilly's type theorem obtained in \cite{BCP}. 

\begin{teor}[cf. Th. 1.5 in \cite{BCP}]
\lb{TF} Let $M$ be an $n$-dimensional compact submanifold of the euclidean space $\re^{n+k}$ endowed with a density $e^\psi$. The first eigenvalue $\lambda_1$ of the $\psi$-Laplacian $\Delta_\psi$ on $M$ associated to the metric and density induced on $M$ by those of $\re^{n+k}$ satisfies the inequality
\begin{equation}\lb{ineqTF}
\lambda_1 \le  \fracc{\int_M |\vHp + \ona\psi|^2 \ \mu_\psi}{n\  V_\psi(M)}.
\end{equation}
Moreover, if equality holds then $M$ is a shrinker for the mean curvature flow in $\re^{n+k}$, and there is a point $p\in \re^{n+k}$ such that $\psi$ restricted to $M$ has the form $\psi|_M = a - \frac12 \lambda_1 r_p^2$, where $r_p$ denotes the Euclidean distance to $p$,  and $\vec{H}(F(x)) = - \lambda _1 (F(x) - p)^{\perp}$, where $F$ is the embedding $F:M\flecha \re^{n+k}$. 

As a consequence, if $k=1$ and $H>0$, or if $n=2$, $k=1$ and $M$ has genus $0$, then equality  holds if and only if $M$ is a round sphere.
\end{teor}

The inequality \eqref{ineqTF} was proved, in a more general form, in Theorem 1.5 in \cite{BCP}. What is new in the statement of the Theorem \ref{TF} is the observation that equality in \eqref{ineqTF} implies that $M$ must be a shrinker of the mean curvature flow in $\re^{n+k}$. This essentially appears when we consider as a condition on the function $\psi|_M$ what in  \cite{BCP} is written as: \lq\lq $M$ is contained in the hypersurface $\lambda \psi+\int^r s_\delta(t) dt =c$" and, moreover, restricts the ambient space to be $\re^{n+k}$.

\begin{nota}  When $\oM=\re^{n+k}$, the inequality of Lemma 3.2 in \cite{BCP} becomes
$\ds
0   = n \ V_{\psi} (M) 
+ \int_M \< \vHp + \overline{\nabla} \psi, F_p \>  \ \mu_{\psi}
$, where $F_p(x) := F(x)-p$.
If $M$ is a shrinker, this formula gives  \ $\ds \frac{\int_M |\vHp + \ona\psi|^2 \ \mu_\psi}{n\  V_\psi(M)} = \fracc{ C^2 \int_M r_p^2 \ \mu_\psi}{n\ V_\psi(M)}$ $=-C$. Thus $-C$  is an eigenvalue of $\Delta_\psi$ on $M$  and the coordinates $x_p^{i} = (x-p)^i$ of $M$  are the corresponding eigenfunctions (see \eqref{xieig}). 

It remains unproved if $-C$ is precisely the first eigenvalue. We shall go back to this question after we state our next theorem.
\end{nota}

\begin{teor}\lb{CorF}
Let $(\oM;\<,\>;\psi)$ be a Riemannian manifold with metric $\<,\>$ and density $e^\psi$. Let $M$ be an $n$-dimensional compact submanifold of $\oM$ of codimension $k$. The first eigenvalue $\lambda_1$ of the $\psi$-Laplacian $\Delta_\psi$ on $M$ associated to the induced metric and density on $M$ satisfies the inequality
\begin{equation}\lb{ineq}
\lambda_1 \le b^2 + \fracc{\int_M |\vHp + \ona\psi|^2 \ \mu_\psi}{n V_\psi(M)} \text{ for some constant $b$ which depends on } \oM
\end{equation}
Moreover, if $\oM$ is a sphere $S^{n+k}$ of radius $R$ in $\re^{n+k+1}$, then $b^2 = n/R^2$ and the  equality implies  $\psi|_M$ is constant, there is a point $p\in \re^{n+k+1}$ such that $\wvH = -\lambda_1 (F-p)^\bot$ (where $\wvH$ is the mean curvature of $M$ in $\re^{n+k+1}$) and, either: 
\begin{description}
\item[ i)] 
 $p$ is the origin $0$ of $\re^{n+k+1}$ and 
$M$ is a minimal submanifold of $S^{n + k}$, or 

\item[ ii)] 
$k=1$ and $M$ is a geodesic sphere of $S^{n + 1}$ with center at $p\ne 0$, or

\item[ iii)] 
$M$ is a  minimal submanifold of  a geodesic sphere $S^{n + k -1}$
of $S^{n + k}$ with center at $p\ne 0$. 
\end{description}
In case (ii) the converse is also true. In cases i) and iii) the converse would be equivalent to Yau's conjecture on the first eigenvalue of a minimal surface of a sphere (see, for instance, \cite{br13} section 5 for the statement of the conjecture and some history). 
\end{teor}

The concept of center of $\psi$-mass, its existence and uniqueness is explained in Section \ref{pre}. For the moment the reader can look at it as a generalization of the center of mass when the volume elements associated to a density are used instead of the standard volume elements.

 Theorem \ref{CorF} has some overlaps with Theorem 1.5 in \cite{BCP}. But there are different results: our inequality \eqref{ineq} is weaker than the corresponding inequality in \cite{BCP}, but valid for 
any submanifold, without restrictions 
on the size of the submanifold. When the ambient space is the sphere, the inequality is the same in both works, but in \cite{BCP} there is restrictions on the size of the submanifold, whereas Theorem \ref{CorF} is valid for any compact submanifold $M$ of $S^{n+k}$. More concisely, Theorems 1.4 and 1.5 in \cite{BCP} generalize Theorems 2.3 and 2.1 in \cite{he}, whereas our theorem \ref{CorF} generalizes Theorem 1.1 in \cite{he}.

In view of the equivalence of the converse in cases i) and iii) of the above theorem with the Yau's conjecture on the first eigenvalue of a minimal surface of a sphere, we conjecture that this is also true for all compact shrinkers in $\re^{n+k}$, that is, we conjecture that if $M$ is a compact shrinker in $\re^{n+k}$ ($\vec{H} =  C \ (F-p)^{\perp}$, $C<0$), then $\lambda_1 = -C$. In \cite{dixi} it is proved that $\lambda_1\in[-C/2,-C]$, and other related bounds for $\lambda_1$ are obtained in \cite{madu} and \cite{piri}. Those estimates are in agreement of our conjecture, but, as a better
support of it, we give the following theorem, which states it in the best known example of shrinker which is not contained in a sphere.

\begin{teor}\lb{l1A} Assume $\psi$ is of Gaussian type, and let $M$ be an embedded shrinker of revolution (around an axis containing $p$) in $\mathbb{R}^3$. That is, $\psi(x)=a+\frac12 C \ r_p(x)^2$ with $a\in \re$ and $C<0$, and $\vec{H_\psi}=0$ on $M$. 

Further assume that the intersection of $M$ with the sphere of radius $\sqrt{2/(-C)}$ centered at $p$ has only two connected components, and that $M$ is symmetric with respect to reflection across the axis of revolution (at $p$). Then the first eigenvalue $\lambda_1$ of the $\psi$-Laplacian $\Delta_\psi$ on $M$ is $\lambda_1 = -C=\ds\fracc{\int_M |\vHp + \ona\psi|^2 \ \mu_\psi}{n\  V_\psi(M)}$. 
\end{teor} 

We have stated the hypotheses of Theorem \ref{l1A} in a general manner, but it is possible that, apart from the sphere, they refers only to a unique (up to homothety) surface - Angenent's torus. In fact, M\o ller \cite{mo11} proves the existence of a torus of revolution satisfying the hypotheses of Theorem \ref{l1A}, and conjectures that this torus is unique (up to a homothety). In fact, it is conjectured that Angenent's torus is the unique embedded shrinker of revolution. The hypothesis on the intersection with the sphere is explicitly proved in Proposition 2.1 (3) of \cite{mo11} for M\o ller's version of the shrinker. In Theorem \ref{an} of this paper we shall prove that Angenent's version of this torus also satisfies this hypothesis. 

What is striking in these results is the fact that equality requires a Gaussian density and a shrinker. In previous theorems on bounds of the $\lambda_1$ for the $\psi$-Laplacian (\cite{huru}) the equality holds only when $\psi$ is constant. 
We found especially interesting the fact that, when $\psi$ is a Gaussian density, for the $\psi$-Laplacian, both, the sphere and the  Angenent's torus  have the same first eigenvalue (Theorem \ref{l1A}).

{\bf Acknowledgments:} 
Research partially supported by  the
DGI (Spain) and FEDER  project MTM2016-77093-P and the Generalitat Valenciana Project  AICO/2019/191.

JZ is partially supported by the NSF under grant number DMS-1802984. 

\section{Preliminaries}\lb{pre}

\begin{defi} Let $f: \Bbb R \flecha \Bbb R$ be a ${C}^2$-function.
Let $M$ be a compact submanifold of $\re^m$ endowed with a density $e^\psi$, and let us denote by $\mathcal{F}: \re^m \flecha \Bbb R$ the function defined by $p \mapsto \mathcal{F}(p) = \int_M f(r_p(x)) \ \mu_{\psi}(x)$. The point  $x_0 \in \Bbb \re^m$ where $\mathcal F$  attains its minimum is called the $f$-center of $\psi$-mass of $M$.
\end{defi}

The following theorem of existence and uniqueness of the $f$-center of $\psi$-mass can be proved following the steps for the usual center of mass (see \cite{BeGo},\cite{Ka}, \cite{BK} and \cite{GMO})

\begin{teor}\lb{com}
Let $f: \Bbb R \flecha \Bbb R$ be a ${C}^2$-function satisfying $f'(t) > 0$ and $f''(t) > 0$ for every $t > 0$. If $M$ is a compact submanifold of $\re^m$, then there is a unique $f$-center of $\psi$-mass of $M$. 
\end{teor}

Let $p$ be the $f$-center of $\psi$-mass of $M$. Because $p$ is a minimum of $\mathcal F$ , the gradient of $\mathcal{F}$ must be $0$ at $p$, that is, for every $\zeta \in \re^m$, $0 = \<\grad\mathcal{F}(p),\zeta\> = \int_M f'(r_p(x)) \< \overline{\nabla} r_p, \zeta  \> \ \mu_{\psi}(x)$. In particular, if we take  $f(t) = \ds\frac{1}{2} t^2$,  $x_p^i$ are  the standard coordinates of $\re^m$ with origin at $p$, and $\zeta = \partial _{i}$ in $\re^m$, we have the following equality 
\begin{equation}\lb{xip0}
0 = \int_M r_p \< \overline{\nabla} r_p, \partial _{i}  \> \ \mu_{\psi} 
= \int_M x^{i}_p  \ \mu_{\psi} ,
\end{equation}
and we can use the $x_p^i$ as test functions to apply the Rayleigh priciple to the $\psi$-Laplacian. From now one we shall refer to the $\ds\frac{1}{2} t^2$-center of $\psi$-mass as the center of $\psi$-mass.

The usual divergence Theorem and the formula of integration by parts  take the following form for $\Delta_\psi$. 
\begin{equation}
\int_{M} {\mathrm{div}}_{\psi} X \  \mu_\psi = 0
\label{intdiv}
\end{equation}

\begin{align}
\int_{M} f {\Delta}_{\psi} h \  \mu_\psi &= 
\int_{M}  h {\Delta}_{\psi} f \  \mu_\psi = - \int_{M}  \< \nabla h, \nabla f \>  \mu_\psi 
\label{intlap}
\end{align}

 Thanks to formula \eqref{intlap}, the basic properties of the eigenvalues and eigenfunctions of $\Delta_\psi$ are the same that for the usual Laplacian (see, for instance, \cite{Ch}, pages 8 and 9), just changing the volume form $\mu$ by $\mu_\psi$. 
 
 On $C^\infty(M)$ we have the scalar product $(f,h)_{\psi 1} = \int_M f h \mu_\psi + \int_M \<\nabla f, \nabla h\> \mu_\psi$. 
 Instead of $L^2(M)$, when there is a density $e^\psi$ we use the space 
$L^2_\psi(M) = \{ f:M\rightarrow \re \ ;\  \int_M f^2 \mu_\psi< \infty\}$.
  One defines also $\mc{H}_\psi(M)$ as the $|| \cdot ||_{\psi 1}$-completion of $\{f\in C^\infty(M);\ ||f||_{\psi 1}<\infty\}$ in $L^2_\psi(M)$.

The standard Rayleigh's principle is still true for $\psi$-Laplacians just changing the volume elements as above.

\section{Proof of Theorem \ref{TF}}\lb{sec12}
 Let $p \in M$ be the center of $\psi$-mass of $M$. Using Rayleigh's principle and the coordinate functions $x^i_p$ as test functions, it is proved in subsection 3.2 of \cite{BCP} that the inequality \eqref{ineqTF} holds. We alert the reader that Batista, Cavalcante and Pyo call $-f$ what we have denoted by $\psi$. Moreover, they observe that the equality in \eqref{ineqTF} implies
  that there is a constant $C$ such that
\begin{equation}\lb{56}
\vHp + \overline{\nabla} \psi = C \ F_p
\end{equation}
The tangent part of the equality \eqref{56} is
\begin{equation}\lb{57}
\nabla \psi = C \ r_p \nabla r_p,
\end{equation}
and its normal part is
\begin{equation}\lb{57-1}
\vHp + \left( \overline{\nabla} \psi \right) ^{\perp}
= C \ (F_p)^{\perp} \qquad \text{ that is } \qquad 
\overrightarrow{H} = C \ (F_p)^{\perp},
\end{equation}
which is the equation of a self-similar solution of the mean curvature flow.

If $C = 0$,  $\overrightarrow{H} = 0$, the submanifold $M$ is minimal and compact in $\re^{n+k}$, which is not possible. 

If $C > 0$,  $M$ is a compact expander. In this case  $M$ will expand  with time under mean curvature flow, then  the volume of $M$ will grow, which is impossible because the mean curvature flow decreases volume. 

If $C < 0$, $M$ is  a shrinker, as claimed in the statement of Theorem \ref{TF}. 
Moreover, from (\ref{57}) 
\begin{align}\lb{57-2}
\nabla \psi = C \ r_p \nabla r_p = \frac{1}{2}\ C \ \nabla r_p^2
\end{align}
then \ $\nabla (\psi - \frac{1}{2}\ C \ r_p^2) = 0$, \ therefore  
\ $\psi|_M = a + \frac{1}{2}\ C \ r_p^2$ for some constant $a$. That is, equality in \eqref{ineqTF} also implies that the density $e^\psi$ restricted to the submanifold $M$ is of Gaussian type.

Now, we check that, under conditions \eqref{57-1} and \eqref{57-2}, the $x_p^i$ are eigenfunctions of $\Delta_\psi$ with eigenvalue  $-C$. In fact
\begin{align}
\Delta _{\psi} x_p^{i} &= \Delta x_p^{i} + \< \nabla \psi, \nabla x_p^{i} \> = \vec{H}^{i} + \< \nabla \psi, \nabla x_p^{i} \>  \nn \\ 
& 
= C \ ((F-p)^\bot)^i    + C \ ((F-p)^{\top})^{i} 
= C \ (F-p)^{i} = C \ x_p^{i}, \lb{xieig}
\end{align}
Then, if we have equality in \eqref{ineqTF}, $M$ is a shrinker soliton ($\vec{H}= C F_p^\bot$, with $C<0$), $\psi$ restricted to $M$ is of Gaussian type, and 
$\lambda _1 = - C$. All these are necessary conditions to have equality in \eqref{ineqTF}, but we have not proved that they are sufficient. We know that in a shrinker $-C$ is an eigenvalue of $\Delta_\psi$ with $\psi$ Gaussian, but we do not know if $-C$ is precisely its first eigenvalue. We shall prove that this is true for some tori of Angenent type in Theorem \ref{l1A} and we conjecture that it is true for all the shrinkers.

When $k = 1$ and $H > 0$, it is known  (see \cite{hui}) that the unique compact shrinking soliton is the sphere, where we know that the equality holds in \eqref{ineqTF}. Also, for $k=1$, $n=2$ and $M$ of genus $0$, it has been proved recently in \cite{br15} that $M$ must be a sphere. Then, when $k=1$ and $H>0$ or $n=2$ and genus of $M$ is zero, the theorem can be stated with a \lq\lq if and only if"\ for the case of equality in \eqref{ineqTF}. 

\bigskip

\section{Proof of Theorem \ref{CorF}}

In this setting $M$ is a $n$-dimensional submanifold of a Riemannian $(n+k)$-dimensional riemannian manifold $\oM$. By the Nash embedding theorem, we can consider $\oM$ as a Riemannian submanifold of $\re^m$. We shall denote by $D$ the covariant derivative and the gradient in $\re^m$. If $\oM$ has a density $e^\psi$, we can consider on $\re^m$ an extension of $\psi$ that we shall denote by  $\wps$. We choose this extension of $\psi$ satisfying 
\begin{equation}\lb{61}
D \wps|_\oM = \overline{\nabla} \psi.
\end{equation}
This can be done, for instance, taking $\wps(p+ t \zeta) = \psi(p)$ for $p+t\zeta$ in a small tubular neighborhood of $\oM$, where $p\in \oM$ and $\zeta$ is a unit vector normal to $\oM$ at $p$. 

We have then the following chain of extrinsic curvatures: 

$\alpha $ \ is the second fundamental form of $M$ in $\overline{M}$ with $\psi$-mean curvature $\vHp$.

$\alpha_{\overline{M}} $ \ is the second fundamental form of  $\overline{M}$ in 
$\re^m$.

$\widetilde{\alpha} $ \ is the second fundamental form of  $M$ in 
$\re^m$, with $\wps$-mean curvature $\wHp$.

Since  $M$ is also a Riemannian submanifold of $\re^m$,  we can apply Theorem \ref{TF} to obtain
\begin{equation}\lb{62}
\lambda_1 (M) \le  \fracc{\int_M |\wHp + D\wps|^2 \ \mu_\psi}{n\  V_\psi(M)}
\end{equation}

Let $\{e_i\}_{i = 1}^n$ be an orthonormal local frame of $M$. From the definitions of the $\alpha$'s and $\vec{H}$'s it follows
\begin{align} \lb{60}
& \wvH = \ds\sum_{i = 1}^{n} \wal(e_i, e_i) 
= \ds\sum_{i = 1}^{n} \alpha_{\overline{M}} (e_i, e_i)   + \vec{H}
\end{align}

From (\ref{61}), (\ref{62}) and \eqref{60}, taking into account that  $\vHp$ and 
$\ona\psi$ are tangent to $\overline{M}$ and 
$\ds\sum_{i = 1}^{n} \alpha_{\overline{M}} (e_i, e_i)$ is orthogonal to 
$\overline{M}$, we obtain  
\begin{align} \lb{61-1}
|\wHp + D \wps|^2  
= |\vHp + \ona\psi|^2 + \left|\ds\sum_{i = 1}^{n} \alpha_{\overline{M}} (e_i, e_i) \right| ^2
\end{align}

Denoting by $c = \ds\max_{v \in TM, |v|=1} \{ \  | \alpha_{\overline{M}}  (v, v) |\  \} $, we have the estimate
\begin{align} \lb{61-2}
\ds\int_M  \left|\ds\sum_{i = 1}^{n} \alpha_{\overline{M}} (e_i, e_i) 
\right| ^2 \ \mu_\psi \leq  n^2 \ c^2 \ V_\psi(M)
\end{align}


By substitution of (\ref{61-1}) in (\ref{62}), having in mind (\ref{61-2}), we obtain 
\begin{align}\lb{63}
\lambda_1 (M) & \le   n \ c^2 
+  \fracc{\int_M |\vHp + \ona\psi|^2 \ \mu_\psi}{n\  V_\psi(M)}  
\end{align}
which is \eqref{ineq}, with $b^2 = n c^2$. 

Now, let us consider that  $\oM$ is a round sphere $S^{n+k}$ of radius $R$. Then we can take $m=n+k+1$ and the sphere centered at the origin of $\re^m$. In this case the inequality \eqref{61-2} is a equality and $c=1/R$.  Let us suppose that, in this setting, we have the equality in \eqref{63}, then we have equality in \eqref{62}, which, by Theorem \ref{TF}, implies that  $M$ is a shrinker in $\re^m$ satisfying $\wvH = -\lambda_1 F_p$ and $\wps|_M = a - \frac12 \lambda_1 r_p^2$ for some point $p\in\re^{n+k+1}$. 

Now we consider two possibilities:

i) $p$ is the center of $S^{n+k}$, then $\wps|_M$ is constant and $M$ is a shrinker of $\re^{n+k+1}$ contained in $S^{n+k}$. But, as Smoczyk noticed (\cite{sm05}), the shrinkers contained in a sphere are the minimal submanifolds of the sphere. Then part i) of Theorem \ref{CorF} is proved.

ii) $p$ is not the center of $S^{n+k}$.
  Then 
\begin{equation} \lb{10*}
\wvH = -\lambda_1 (F_p)^{\perp} = -\lambda_1 (F - p)^{\perp} = -\lambda_1 \< F - p, N \> N -\lambda_1  \< F - p, \zeta \> \zeta   
\end{equation}  
where $\perp$ means the component orthogonal to $M$ in $\Bbb R^{ n + k +1}$, $N$ is the inner unit vector orthogonal to $S^{n+k}$ in $\re^{n+k+1}$ and $\zeta$ is a unit vector in the direction of the component of $F-p$ orthogonal to $M$ in $S^{n+k}$.

In the sphere, \eqref{60} becomes 
\begin{align}\lb{11*}
\wvH =  \frac{n}{R} \ N + \vec{H}, 
\end{align}
From (\ref{10*}) and (\ref{11*}) we obtain 
\begin{align}\lb{12*}
&\frac{n}{R} =   - \lambda_1 \< F - p, N \>
\end{align}
Now we show that 
(\ref{12*}) implies that $M$ is contained in a sphere 
of dimension $n+k-1$ contained in $\oM=S^{n+k}$. In fact, let $\beta(s)$ be a curve in $M$, by (\ref{12*}) one has $\ds\frac{- n}{ \lambda_1 R} = \< \beta(s) - p, N \>$.
Taking derivative respect to $s$, and using that $\beta(s)$ is also contained in $S^{n+k}$,
\begin{equation}
 0 = \< \beta', N \> + \< \beta - p, (N \circ \beta)' \> = - \< \beta - p, \frac{1}{R} \beta' \>. 
\end{equation}
Then  $0=\< \beta - p, \beta' \> =  \ds\frac{1}{2} \< \beta- p, \beta - p \> '= \fracc12 \ds\frac{d}{ds}\left( |\beta - p|^2 \right)$, so $|\beta - p|$  is constant, and $\beta(s)$, therefore $M$,  is contained in the intersection of a sphere $\mathbb S$ of dimension $n+k$ centered at   $p$  with $S^{n+k}$, which 
 is a sphere $S^{n+k-1}$ of dimension $n+k-1$.
If $k=1$, $M$ is just the sphere where it is contained. If $k>1$,   since $M$ is a shrinker contained in a sphere $\mathbb S$ centered at $p$, $M$ is minimal in $\mathbb S$ and it is contained in the sphere $S^{n+k-1}\subset \mathbb S$, but this is possible only if $S^{n+k-1}$ is totally geodesic in $\mathbb S$, which is equivalent to say that $p$ is also the center of $S^{n+k-1}$ and $M$ is minimal in $S^{n+k-1}$. 
This finishes the proof of parts ii) and iii) of Theorem \ref{CorF}.

\bigskip

\section{Proof of Theorem \ref{l1A}}
In this section we shall work in $\re^3$. Then, for simplicity, we shall use $x$, $y$, $z$ to denote the coordinates $x^1$, $x^2$, $x^3$.

It was proven by Brendle in \cite{br15} that the unique embedded shrinker of genus $0$ in $\re^3$ is the standard round sphere. Since an embedded surface of revolution only can be of genus $0$ and $1$, and Theorem \ref{l1A} is well known for the round sphere, we shall suppose that the embedded shrinkers in the hypothesis of the theorem are tori of  revolution $T$ around the $z$-axis. They are warped products  $T = E \times_{\rho} S^1$ with $E$ a simple closed curve in the $xz$-plane, invariant respect to the symmetry   $z \longrightarrow - z$, where $\rho: E \flecha \Bbb R$ denotes the distance, in the plane, to the $z$-axis, and also the distance $\rho_3$ in $\re^3$ to the $z$-axis, restricted to $T$. From now on we shall use the same letter for $\rho$ and $\rho_3$ restricted to $E$ and on $T$ respectively, because $\rho_3(s,\theta) = \rho(s)$ for every $(s,\theta)\in E\times_\rho S^1$. 

Using the expression of the standard Laplacian on a warped product (see, for instance, \cite{Ej}), the $\psi$-Laplacian on the torus $T$ can be written as 
\begin{align}\label{lap}
\Delta_{\psi} ^T &= \Delta ^{E} +  \frac{1}{\rho} \ \nabla{\rho}
+  \frac{1}{\rho^2} \ \Delta ^{S^1} +   \nabla \psi \\
& = \Delta _{\ln \rho + \psi}^{E} + \frac{1}{\rho ^2} \Delta ^{S^1}, \label{lap2}
\end{align}
where $\Delta^{E}$ and $\Delta^{S^1}$ denote the usual Laplacians in $E$ and $S^1$ respectively. 

To study the spectrum of $\Delta_\psi^T$  we shall follow the procedure of Ejiri \cite{Ej} on warped products.
 Let $\{ \psi_0, \psi_1, \psi_2,  \dots\  \}$ be a complete orthonormal  basis of $L^2(S^1)$ formed by eigenfunctions of $\Delta^{S^1}$ with corresponding eigenvalues  $\lambda_0=0<\lambda_1 =\lambda_2 =1 <\lambda_3 =\lambda_4 =4 < \cdots< \lambda_{2k-1}=\lambda_{2k}=k^2 < \cdots $.  
For every $\lambda_i$, we define the operator  $L_{\lambda_i}^{\psi}$ acting on smooth real functions defined on $E$
 \begin{align}\lb{llp}
 L_{\lambda_i}^{\psi} = \Delta ^{E} +  \ds\frac{1}{\rho} \ \nabla \rho  - \ds\frac{\lambda_i}{\rho^2} +  \nabla \psi .
 \end{align}
The same arguments given in Lemma 2.3 of \cite{Ej}  show that $L_{\lambda_i}^{\psi}$ is a strongly elliptic self-adjoint operator on the space $L_{\ln\rho+\psi}^2 (E)$  of the  smooth functions on $E$ with bounded square norm under the measure $\rho e^\psi \mu = dv_{\ln\rho+\psi}$. Then, there exists a complete orthonormal basis of $L_{\ln\rho+\psi}^2 (E)$ formed by eigenfunctions $\phi_j^{\lambda_i}$ of $L_{\lambda_i}^{\psi}$ with eigenvalues  $\mu_0^{\lambda_i} \le \mu_1^{\lambda_i}\le  \mu_2^{\lambda_i}  \le \cdots $.

Now, we use \eqref{lap} to compute
\begin{align}\lb{dphit}
\Delta_{\psi} ^T \left( \psi_i \phi_j^{\lambda_i} \right)   \nonumber 
& = \left( \Delta ^{E} +  \frac{1}{\rho} \ \nabla \rho  - \frac{\lambda_i}{\rho^2} + \frac{\lambda_i}{\rho^2} 
+  \frac{1}{\rho^2} \ \Delta ^{S^1} +  \< \nabla \psi, \nabla \cdot \>\right) \left( \psi_i \phi_j^{\lambda_i} \right)  \\ 
& = L_{\lambda_i}^{\psi} \left( \psi_i \phi_j^{\lambda_i} \right) 
+ \left( \frac{\lambda_i}{\rho^2} +  \frac{1}{\rho^2} \ \Delta ^{S^1} \right) 
\left( \psi_i \phi_j^{\lambda_i} \right) = L_{\lambda_i}^{\psi} \left( \psi_i \phi_j^{\lambda_i} \right) 
\end{align}
because $\phi_j^{\lambda_i}$ depends only on  $E$, $\psi_i$ only depends on $S^1$ and
$\Delta ^{S^1} \psi_i = - \lambda _i \psi_i$. 
Moreover, 
\begin{align}\label{pg3 (22-9)}
L_{\lambda_i}^{\psi} \left( \psi_i \phi_j^{\lambda_i} \right) 
= \Delta ^{E} \left( \phi_j^{\lambda_i}  \right) \ \psi_i + \frac{1}{\rho} \<\nabla \rho , \nabla\phi_j^{\lambda_i}  \> \ \psi_i - \frac{\lambda_i}{\rho^2} \ \psi_i \phi_j^{\lambda_i} + \< \nabla \psi, \nabla \left(\psi_i \phi_j^{\lambda_i}\right) \>
\end{align}

But $\nabla \psi$ is in the direction of the position vector, this is orthogonal to the curves $S^1$ in the warped product, and $\nabla \psi_i$ is tangent to these $S^1$, then $\<\nabla \psi, \nabla\psi_i\>=0$ and 
$\< \nabla \psi, \nabla \left(\psi_i \phi_j^{\lambda_i}\right) \> = \< \nabla \psi, \nabla (\phi_j^{\lambda_i}) \> \psi_i $.

Plugging this equality into \eqref{dphit} and \eqref{pg3 (22-9)} and taking into account that  $\mu_j ^{\lambda_i}$ is an eigenvalue of $L_{\lambda_i}^{\psi}$ with eigenfunction $\phi_j^{\lambda_i}$ 
\begin{align}\label{eigT}
\Delta_{\psi} ^T \left( \psi_i \phi_j^{\lambda_i} \right)= L_{\lambda_i}^{\psi} \left( \psi_i \phi_j^{\lambda_i} \right) 
& = - \mu_j ^{\lambda_i} \phi_j^{\lambda_i} \psi_i
\end{align}
that is, $\psi_i \phi_j^{\lambda_i} $ are eigenfunctions of $\Delta _{\psi} ^T$ with eigenvalues $\mu_j ^{\lambda_i}$. 
From this, using the same arguments that in the proof of Theorem 2.1 in \cite{Ej},  we conclude that $\{\psi_i \phi_j^{\lambda_i} ; \ i,j=0,1,2,\dots\}$  is a complete orthonormal basis of $L^2_\psi(T)$ and the $\{\mu_j^{\lambda_i}\}$ are all the eigenvalues of $\Delta_\psi^T$.

To finish the proof of Theorem \ref{l1A} it remains only to check  that $\inf\{\mu_j^{\lambda_i}\ne 0; \ i,j=0,1,2,\dots\} = - C$. We shall prove this in a series of lemmas.

\begin{lema}\lb{e0Llp} Constant functions are eigenfunctions of $L_{\lambda_i}^\psi$ with eigenvalue $0$ if and only if $i=0$
\end{lema}
\begin{demo}
If we write $L_{\lambda_i}^\psi$ in the form
\bec\lb{LlnD}
L_{\lambda_i}^\psi =  \Delta _{\ln \rho + \psi} ^{E} - \frac{\lambda_i}{\rho ^2},
\eec
by \eqref{LlnD} and \eqref{intlap}, we have 
\begin{align}\label{pg.8 (29-9)}
 \int_E \(L_{\lambda_i}^\psi u\) \ u \ dv_{\ln \rho + \psi} 
&= \int_E \left( \Delta _{\ln \rho + \psi} ^{E} u \right)  \ u \ dv_{\ln \rho + \psi} 
- \int_E  \frac{\lambda_i}{\rho ^2} \ u^2 \ dv_{\ln \rho + \psi} \nn\\ &= - \int_E |du|^2 \ dv_{\ln \rho + \psi} 
- \int_E  \frac{\lambda_i}{\rho ^2} \ u^2 \ dv_{\ln \rho + \psi} \le 0
\end{align}
which vanishes if and only if $u$ is constant and $i=0$. 
\end{demo}

In the next two lemmas we check, adapting the exposition given in \cite{Ch},  that the well known Nodal Domains Theorem of Courant is still valid for a Laplacian with density. First we recall the concept of nodal domain and the classical Sturm theorem.

The {\it nodal domains of a function} $f:M\flecha \re$ are the connected components of the set $M\setminus \{ x \in M  \, |  \ f(x) = 0 \}$.

\textbf{Sturm's Theorem} (\cite{CL}, page 214) {\it The eigenvalue problem $(p y')'+ \mu \ p\ y =0$ for periodic functions $y$ on $\re$ of period $L$  has solution for a sequence of values of $\mu$
$$ 0 = \mu_0 < \mu_1 \le \mu_2 < \mu_3 \le \mu_4 < \cdots < \mu_{2j-1} \le \mu_{2j} < \cdots$$
where $\mu_0$ has multiplicity $1$, if $\mu_{2j-1} < \mu_{2j}$ both $\mu_{2j-1}$ and  $\mu_{2j}$ have multiplicity $1$, if $\mu_{2j-1} = \mu_{2j}$, this eigenvalue has multiplicity $2$. Moreover, the eigenfunctions corresponding to $\mu_{2j-1}$ and $\mu_{2j}$, for $j\ge 1$,  have $2 j$ nodal domains.} 

In Lemma \ref{lema1} we shall write the above eigenvalue problem as the eigenvalue problem of the operator $y\mapsto \fracc{(p\ y')'}{p}$. We shall call it the {\it Sturm operator}.

\begin{lema}\lb{Cour}
 Let $\{\phi_0, \phi_1, \phi_2, \cdots \}$ be a complete orthonormal basis of $L^2_\varphi(M)$ with each $\phi_j$ an eigenfunction of $\Delta_\p$ associated to an eigenvalue $\lambda_j$, $j=0,1,2, \cdots $. Then the number of nodal domains of $\phi_k$ is less or equal to $k+1$, for every $k=0,1,2,\cdots$.  Moreover, the number of nodal domains of $\phi_0$ and  $\phi_1$ are $1$ and  $2$ respectively.
\end{lema}
\begin{demo}
It follows with exactly the same arguments used for the ordinary Laplacian (see \cite{Ch} pages 19-20), just changing the corresponding spaces associated to the Riemannian measure $\mu_g$ by those spaced associated to the measure with density $\mu_\psi$ indicated in the Preliminaries.
\end{demo}

With these lemmas in hand, we can give the final steps in the proof of the theorem by comparing to some known eigenfunctions.

\begin{lema}\lb{lema1} $\inf\{\mu_j^{\lambda_i}\ne 0; \ i,j=0,1,2,\dots\} = \min\{\mu_1^{\lambda_0},\  \mu_2^{\lambda_0},\ \mu_0^{\lambda_1}\}$. 

Moreover,  $\mu_0^{\lambda_1}=-C$.
\end{lema}
\begin{demo}
By Lemma \ref{Cour}, the eigenfunctions associated to the first nonzero eigenvalue of $\Delta_\psi^T$ have  exactly two nodal domains. Then we shall look for  products $\psi_i\phi_j^{\lambda_i}$ with two nodal domains.

Let us denote by $N(f)$ the number of nodal domains of a function. Obviously we have $N(\psi_i\phi_j^{\lambda_i})= N(\psi_i) N(\phi_j^{\lambda_i})$. We discuss the possibilities for the different values of $i$ and $j$.

\begin{enumerate}
\item For $i=0$,  $\psi_0$ is constant, then $N(\psi_0 \phi_j^{\lambda_0})= N(\phi_j^{\lambda_0})$. Moreover, from \eqref{llp}, $L_{\lambda_0}^{\psi} (y) = (\Delta ^{E} +  \ds\frac{1}{\rho} \ \nabla \rho   +  \nabla \psi) y= y'' + \(\frac{\rho'}{\rho}+\psi'\)y'$ has the form of  the Sturm operator, with $p = \rho e^\psi$. Then, by Sturm's Theorem, the eigenfunctions  $\phi_{2j-1}^{\lambda_0}$ and $\phi_{2j}^{\lambda_0}$ have $2j$ nodal domains for $j\ge 1$, whereas, by Lemma \ref{e0Llp},  $\phi_{0}^{\lambda_0}$ is constant and has one nodal domain.  By Lemma \ref{Cour} the first nonzero eigenvalue of $\Delta^T_\psi$ has an associated eigenfunction with two nodal domains. Therefore the unique candidates to be this eigenvalue when $i=0$ are $\mu^{\lambda_0}_1$ and $\mu^{\lambda_0}_2$.
\item For $i\in\{1,2\}$, $N(\psi_1)=2$, then $\mu^{\lambda_1}_j=\mu^{\lambda_2}_j$ is a candidate to be the first nonzero eigenvalue of $\Delta_\psi^T$ only if  $N(\phi_{j}^{\lambda_1})=1$. 
Since $T$ a torus of revolution around the $z$-axis and is symmetric respect to reflection across the $xy$-plane,  $T$  has central symmetry and, because $\psi$ is radial, the origin is the center of $\psi$-mass of $T$. It follows from  \eqref{xieig} that the coordinate functions $x,y,z$ restricted to $T$ are eigenfunctions of $\Delta_\psi^T$ with eigenvalue $-C$. Considering the relation $x(s,\theta) = \rho(s) \cos \theta , \quad 
y(s,\theta) = \rho(s) \sin \theta,$  $s\in E$  and $\theta\in S^1$   in the surface of revolution $T$, using  the expression \eqref{lap2} and the fact  that $\cos \theta$ is an eigenfunction of $\Delta^{S^1}$ with eigenvalue $\lambda_1=1$, we obtain 
\begin{align}
C \rho(s) \cos \theta &= \Delta_{\psi}^T (\rho \ \cos \theta) = \left( \Delta _{\ln \rho + \psi}^{E} + \frac{1}{\rho ^2} \Delta ^{S^1} \right) (\rho \ \cos \theta)\nn \\
& = \cos \theta \ \Delta _{\ln \rho + \psi}^{E} \rho  + \frac{\rho}{\rho ^2} \Delta ^{S^1} \cos \theta = \(\Delta _{\ln \rho + \psi}^{E} \rho - \frac1\rho\) \cos \theta . \lb{Cc}
\end{align}

It follows from \eqref{Cc} that 
\begin{align}\label{33}
L_{\lambda_1}^{\psi} \rho =  \Delta _{\ln \rho + \psi}^{E} \rho - \frac{1}{\rho}  =   C \ \rho ,
\end{align}
that is $\rho$ is an eigenfunction of $L_{\lambda_1}^{\psi}$. Moreover $\rho>0$ because the profile curve $E$ does not touch the $z$-axis, hence $N(\rho)=1$, and any other eigenfunction of $L_{\lambda_1}^{\psi}$ orthogonal to $\rho$ must have at least $2$ nodal domains. We conclude that the unique $\mu_j^{\lambda_1}$ candidate to be the first nonzero eigenvalue of $\Delta_\psi^T$ is $\mu_0^{\lambda_1}=-C$ (with $j=0$ because these eigenvalues are ordered), and in fact the eigenfunction $\phi_0^{\lambda_1}$ is a multiple of $\rho$.
\item For $i\ge 3$, $N(\psi_i)\ge 4$, then none of the $\mu_j^{\lambda_i}$ is a candidate to be the first nonzero eigenvalue of $\Delta_\psi^T$.
\end{enumerate}
 \end{demo}

\begin{lema} $\mu_1^{\lambda_0} = -C < \mu_2^{\lambda_0}  = - 2 C$.
\end{lema}

\begin{demo}
As we observed in the previous lemma, it follows from \eqref{xieig} that  $z$ is an eigenfunction of $\Delta^T_\psi$ with eigenvalue $-C$.  This, together with  the expressions \eqref{lap} and \eqref{llp} and the facts that $\lambda_0=0$ and $z$ restricted to $S^1$ is constant give that $z$ is an eigenfunction of $L_{\lambda_0}^\psi$ with eigenvalue $-C$. Since the profile curve $E$ of $T$ is symmetric respect to the $x$-axis, the function $z$ restricted to $E$ has two nodal domains. 

On the other hand, we can verify that $\rho^2 + z^2 + \frac{2}{C}$ is an eigenfunction of $\Delta^T_\psi$ with eigenvalue $- 2C$. In fact, $r_p^2 = \rho^2 + z^2$ and 
\begin{equation*}
\Delta^T_\psi r_p^2 = \sum_{i = 1}^{2} \Delta^T_\psi (x^{i})^2 
= \sum_{i = 1}^{2} 2 x^{i}  \Delta^T_\psi x^{i} + \sum_{i = 1}^{2} 2 | \nabla x^{i}|^2
= 2 C \sum_{i = 1}^{2} (x^{i})^2 + 2 \sum_{i = 1}^{2}  | \nabla x^{i}|^2 
= 2 C r_p^2 + 4
\end{equation*}
then 
\begin{equation*}
\Delta^T_\psi \left( r_p^2 + \frac{2}{C} \right) = 2 C \left( r_p^2 + \frac{2}{C} \right) .
\end{equation*}

But by the hypothesis that the intersection with the sphere of radius $\sqrt{\frac{2}{-C}}$ has two connected components, the function 
$\rho ^2 + z^2 + \frac{2}{C}$ restricted to $E$ also has two nodal domains. Then, by Sturm's Theorem (see Case 1 in the proof of Lemma \ref{lema1}), $\mu_1^{\lambda_0} = -C $ and 
$ \mu_2^{\lambda_0}  = - 2 C$. The corresponding eigenfunctions $\phi_1^{\lambda_0}$ and $\phi_2^{\lambda_0}$ are proportional to $z$ and $\rho^2+z^2+\frac{2}{C}$ respectively. 
\end{demo} 

With this lemma finish the proof of Theorem 3. 

\

Angenent's original construction of a shrinker torus in \cite{An} is obtained by revolution of a closed curve $\Gamma$
in $\Bbb R^2$, contained in the half space $\rho > 0$. We shall prove

\begin{teor}\label{an}
The circle of radius $\sqrt{\frac{2}{-C}}$ cuts $\Gamma$ exactly twice. 
\end{teor}

\begin{demo}
First we recall Angenent's construction. The curve $\Gamma$ parametrized as $\gamma(t)=(\rho(t),z(t))$ is the solution of differential equation 
\begin{equation}\label{dk2}
\frac{z' \rho'' - \rho ' z''}{(z')^2 + (\rho ')^2} - \frac{1}{\rho} z' =  C (\rho z' - z  \rho ' )
\end{equation}
Moreover there exist $R_*  > \sqrt{\frac{2}{-C}}$ such that for all $r \geq R_*$ the solution of \eqref{dk2}  with the initial conditions $\gamma (0) = r$ and $\gamma ' (0) = (0, \mu)$, $\mu \in \Bbb R, \ \mu > 0$ re-intersects the $\rho$-axis before having tangent vector again parallel to the $z$-axis. When $\gamma (0) = R_*$, $\gamma$ re-intersects the $\rho$-axis in $\rho _* < \sqrt{\frac{1}{-C}}$ at a point $\gamma (t_0)$ with $\gamma '(t_0) = (0, \nu)$, $\nu \in \Bbb R, \ \nu < 0$.

Furthermore, Angenent proves that $\Gamma $ can be described by the graph of a function $z = f(\rho)$ with $\rho \in [\rho_*, R_* ]$
and this reflection over the $\rho$-axis, and  from (\ref{dk2}) $f$ satisfies the differential equation 
\begin{equation}\label{ff}
\frac{f''}{1 + (f')^2} = \left( - C \rho - \frac{1}{\rho} \right) f' + C f
\end{equation}

{\bf Step 1.} {\it If $f(\sqrt{\frac{1}{-C}}) < 0$ and $f' (\rho) \geq 0$ whenever $\rho \geq \sqrt{\frac{1}{-C}}$, then $f(\rho) > 0$ whenever $\rho > \rho_0$; $\rho_0$ being the zero of the function 
\begin{align*}
\xi (\rho) =&   \rho \Ei\left(\frac{-C}{2}\rho^2\right) - \sqrt{\frac{1}{-C}}\Ei\left(\frac{1}{2}\right)  - \sqrt{\frac{2\pi}{-C}} \left( \erfi\left(\sqrt{\frac{-C}{2}}\rho\right) - \erfi\left(\frac{1}{\sqrt{2}}\right) \right)\\&+\Ei\left(\frac{1}{2}\right) \left(\rho-\sqrt{\frac{1}{-C}}\right)- \frac{2 \sqrt{e}}{\sqrt{-C}} 
\end{align*} 
for $\rho > \sqrt{\frac{1}{-C}}$,   where $\Ei (x) = \int_{- \infty} ^x  \frac{e^t}{t} \ dt$ and $\erfi(x) = \frac{2}{\sqrt{\pi}}\int_0^x e^{t^2}dt$. }

When $C = - \frac{1}{2}$, $\rho_0 \approx 3.1769$. 

We prove this step following the ideas in Drugan-Kleene \cite{DK}:

Differentiating \eqref{ff} respect to $\rho$, 

\begin{equation}\label{ff2}
\frac{f'''}{1 + (f')^2} = \left( \frac{2 f' f''}{(1 + (f')^2)^2} - C \rho - \frac{1}{\rho}\right) f'' 
+ \frac{1}{\rho ^2} f'
\end{equation}
If $ - f(\sqrt{\frac{1}{-C}}) > 0$, then for $\rho \geq \sqrt{\frac{1}{-C}}$, $f'(\rho)\geq 0$, we have 
$- C \rho - \frac{1}{\rho} \geq  0$, and hence $f''(\rho) > 0$ and $f'''(\rho) \geq 0$. Moreover, from \eqref{ff2} and $f'(\rho) \geq 0$, 
\begin{equation*}
\frac{f'''}{f''} \geq  \frac{2 f' f''}{1 + (f')^2} - C \rho - \frac{1}{\rho}
\end{equation*}
Integrating we find that 

\begin{align*}
\ln f'' (\rho) - \ln f''\left(\sqrt{\frac{1}{-C}}\right)  \geq & \ln (1 + (f')^2)(\rho) - \ln (1 + (f')^2)\left(\sqrt{\frac{1}{-C}}\right) 
  \\ & - \ln (\rho) + \ln (\sqrt{\frac{1}{-C}}) - \frac{C}{2} \rho^2 
+ \frac{C}{2}\frac{1}{-C}   
\end{align*}
\begin{equation*}
f'' (\rho) \geq f''\left(\sqrt{\frac{1}{-C}}\right) \frac{\frac{1}{\sqrt{-C}} (1 + (f')^2(\rho))}{\rho \left(1 + (f')^2\left(\sqrt{\frac{1}{-C}}\right)\right)} e^{-\frac{1}{2} - \frac{C}{2} \rho^2}
\end{equation*}
but, evaluating \eqref{ff} at $\sqrt{1/-C}$,
\begin{equation*}
\frac{f''\left(\sqrt{\frac{1}{-C}}\right)}{1 + (f')^2\left(\sqrt{\frac{1}{-C}}\right)}  = Cf(\sqrt{\frac{1}{-C}})
 \geq 0,  
\end{equation*}
then 
\begin{equation}\label{ff3}
f'' (\rho) \geq \frac{- \sqrt{-C} f(\sqrt{\frac{1}{-C}})}{\rho} e^{-\frac{1}{2}} e^{- \frac{C}{2} \rho^2}
\end{equation}
Integrating twice we find that

\begin{equation*}
f(\rho) -  f\left(\sqrt{\frac{1}{-C}}\right)  -  f'\left(\sqrt{\frac{1}{-C}}\right) \left( \rho - \sqrt{\frac{1}{-C}} \right)  \geq \frac{- \sqrt{-C}}{2\sqrt{e}}  f\left(\sqrt{\frac{1}{-C}}\right) \alpha(\rho),
\end{equation*}
where 

\begin{align*}
\alpha(\rho)=&
 \rho \Ei\left(\frac{-C}{2}\rho^2\right) - \sqrt{\frac{1}{-C}}\Ei\left(\frac{1}{2}\right)  - \sqrt{\frac{2\pi}{-C}} \left( \erfi\left(\sqrt{\frac{-C}{2}}\rho\right) - \erfi\left(\frac{1}{\sqrt{2}}\right) \right)\\&+\Ei\left(\frac{1}{2}\right) \left(\rho-\sqrt{\frac{1}{-C}}\right).
 \end{align*}

Then because $f'\left(\sqrt{\frac{1}{-C}}\right) \rho - f'\left(\sqrt{\frac{1}{-C}}\right) \sqrt{\frac{1}{-C}} = f'\left(\sqrt{\frac{1}{-C}}\right) \left( \rho - \sqrt{\frac{1}{-C}} \right)  \geq 0$, 
we have 
\begin{equation}\label{ff4}
f(\rho) \geq - f(\sqrt{\frac{1}{-C}}) \left( -1 + \frac{ \sqrt{-C}}{2\sqrt{e}}\alpha(\rho)\right),
\end{equation}

Thus we conclude that $f(\rho) \geq 0$ if and only if   $\xi(\rho) \geq 0$. 

For $\rho >\sqrt{\frac{1}{-C}}$, the function $\xi$ is increasing with $\rho$, hence $\xi(\rho) > 0$ if and  only if $\rho > \rho_0$.

When $C = - \frac{1}{2}$ (which is the value taken by Angenent, Drugan, Moller and Kleene) a simple computation with Mathematica gives $ \rho_0 \approx 3.1769$. 

{\bf  Step 2.}  
{\it If $R_* > \rho _0$, the circle of radius $\sqrt{2/(-C)}$ cuts $\Gamma$ exactly twice. }

Whenever we can write the curve $\Gamma$ as a graph of $\rho$ over $z$, $\rho=u(z)$, it follows from (\ref{dk2})  that it must  satisfy the differential equation
\begin{align}\lb{shrrho}
\frac{u''}{1+u'^2} = -C\ z\ u' + C\ u + \frac{1}{u}
\end{align}
By \cite{KM} (Proposition 2 and Lemma 9) $\Gamma$ has exactly two points with tangent line parallel to the $\rho$-axis, then $\Gamma = \text{graph}(u_1)\cup \text{graph}(u_2)$, where $u_1$, $u_2$ are functions defined on some interval $[-\beta, \beta]$, satisfying $u_1>u_2$  and  \eqref{shrrho} on $]-\beta,\beta[$. 

By \cite{DK} (Proposition 3) the second derivative of each segment $u_i$ changes of sign at most twice, and only changes twice if $u_i$ has two minima. Since $\Gamma$ has tangent parallel to the $z$-axis only when $z=0$, the latter cannot occur. By the symmetry of $u_i$ respect to the $\rho$-axis, if the sign of $u_i''$ changes exactly once, then the change has to be at $z=0$, where also $u_i'(0)=0$. If so, from \eqref{shrrho}, $u_i(0) = \sqrt{-1/C}$. By the uniqueness of the solution of \eqref{shrrho}, it must be $u_i=\sqrt{-1/C}$, which is absurd. Then $u_1$ is a concave function and $u_2$ is convex. This implies that $\text{graph}(u_2)$ intersects the circle $z^2 + \rho^2 = 1/C^2$ at most once in the quadrant $z>0$.

We now claim that  $\text{graph}(u_2)$ intersects the line $\rho = \sqrt{-1/C}$. 
Suppose not. Then, since $\rho_*<\sqrt{-1/C} < R_*$, 
the line $\rho=\sqrt{-1/C}$ has to intersect $\text{graph}(u_1)$. By the concavity of $u_1$, $u_1'(z)\ge 0$ for $z\in]-\beta,0[$ and $\text{graph}(u_1)$ intersects $\rho = \sqrt{-1/C}$ exactly once for $z\in[-\beta,0]$. Switching to the representation of $\Gamma$ as a graph $z=f(\rho)$ used in step 1, and taking a reflection over the $\rho$-axis, this means that $f(\sqrt{-1/C})<0$ and $f'(\rho)\ge 0$ for $\rho\ge \sqrt{-1/C}$. By Step 1, this implies that $f(\rho)>0$ for $\rho >\rho_0$, which is equivalent to say that $\text{graph}(u_1)$ crosses the $\rho$-axis at $R_*\le \rho_0$, in contradiction with our hypothesis in this step.

From the above claim, it follows that if $\text{graph}(u_1)$ crosses the circle $z^2 + \rho^2 = 2/(-C)$, it must cross at $\rho > \sqrt{-1/C}$. The tangent line  to the circle at the point $(\rho=\sqrt{-1/C},z=\sqrt{(-1/C)})$ cuts the $\rho$-axis at $\rho = 2/\sqrt{-C}$. Therefore if $\text{graph}(u_1)$ intersects the circle twice in the quadrant $\rho>0,\ z>0$, by concavity one must have $R_* \le  2/\sqrt{-C}$, which again contradicts the hypothesis of this step. Then  $\text{graph}(u_1)$ crosses the circle $z^2 + \rho^2 = 2/(-C)$ at most once in the quadrant $\rho>0,\ z>0$. Moreover, if $\text{graph}(u_1)$ does cross, then by convexity, $\text{graph}(u_2)$ cannot cross the circle. 

Then, by symmetry, we know that the circle centered at the origin and with radius $\sqrt{2/(-C)}$ can cross only $\text{graph}(u_1)$ or $\text{graph}(u_2)$ and, if crosses one of them in the quadrant $\rho>0,\ z>0$, it crosses only once (that is, at most twice total). On the other hand, $R_*>\sqrt{-2/C}>\sqrt{-1/C}> \rho_*$,  so the circle must cross $\Gamma$, and it crosses exactly twice.

{\bf Step 3}. {\it For $C = - \frac{1}{2}$, \ $R_* > \rho _0 \approx 3.1769$.} 

When $C=-1/2$ we can show by rigorous numerical computation that the solution $(\rho(t), z(t))$ of (\ref{dk2}) for $(\rho (0), z(0)) = (3.3, 0)$ and $(\rho' (0), z'(0)) = (0, 1)$ has tangent vector parallel to the $z$-axis before it intersects the $\rho$-axis. The details of this computation and its analysis are deferred to the appendix.

By the property characterizing $R_*$ stated at the begining of the proof, this implies that $R_* > 3.3 > \rho_0$. 

This proves the theorem for $C = - \frac{1}{2}$. For  other values of $C$ the corresponding shrinkers are  homothetic, so the statement of the theorems are true for any $C$.
\end{demo}

\section*{Appendix}

In this appendix we detail the numerical computation and analysis near Angenent's torus. For $C=-\frac{1}{2}$, the solutions of \eqref{dk2} can be viewed as the solutions to the following first-order ODE system (see also \cite{DK}, eq. (9)):

\begin{eqnarray}
\alpha'(s) &=& \frac{z(s)}{2} \sin \alpha(s) + \left(\frac{1}{\rho(s)} - \frac{\rho(s)}{2}\right) \cos \alpha(s) \\ 
z'(s)&=& \cos\alpha(s)\\
\rho'(s) &=& \sin\alpha(s).
\end{eqnarray}

The goal of this appendix is to prove the following:

\begin{prop}
\label{prop:numerical}
Consider the solution $(\rho(s),z(s),\alpha(s))$, $s\geq 0$, of the above ODE system with initial conditions $(\rho(0),z(0),\alpha(0))= (3.3,0, 0)$. 

Let $s_2$ be the least $s>0$ such that $z(s)\leq 0$. Then there is $s_1<s_2$ for which $\alpha(s_1)<-\pi$. That is, by continuity, the solution curve must become parallel to the $z$-axis before it crosses the $\rho$-axis again. 
\end{prop}

We consider the solution above for $s\in [0,\overline{s}]$, where $\overline{s}=3.78$, 
and note that the initial conditions lie in the region
\begin{equation}\label{eq:apbounds}
\mathcal{R} = \{(\rho,z,\alpha) | \rho \in [0.42,4] , |z| \leq 1 ,\ |\alpha|\leq \pi+1\}. 
\end{equation}
We run a numerical simulation, with theoretically bounded error for solutions lying in the region $\mathcal{R}$; this bound will inductively ensure that our solution indeed stays in the region $\mathcal{R}$ and will also yield the final result. 

\subsection*{RK method and error analysis}
\label{sec:error}

Specifically, we use a fixed-step Runge-Kutta method, and follow the method described in \cite{HNW} to estimate its global error. In principle any method may be used, but we have used the standard third-order (RK3) method which converges in a reasonable time and has sufficiently simple error analysis. 

For an ODE system $Y' = f(Y)$, the RK3 method with step size $h$ proceeds as \[Y_{n+1} -Y_n = \frac{h}{6}(k_1 + 4k_2 + k_3),\] where $k_1 = f(Y_n)$, $k_2 = f(Y_n + hk_1/2)$ and $k_3 = f(Y_n + h(2k_2-k_1))$. 

The local (one-step) error bounded by $|Y_{n+1}-Y_n|\leq Ch^4$, where $C$ is an upper bound for \[\frac{1}{4!}\|Y^{(4)}\|_\infty + \frac{1}{3!}\left( \frac{1}{6} \|k_1^{(3)}\|_\infty+ \frac{2}{3} \|k_2^{(3)}\|_\infty  + \frac{1}{6} \|k_3^{(3)}\|_\infty\right).\] Expanding the derivatives $k_i^{(3)}$ with the chain rule, we have $k_1^{(3)}=0$, $\|k_2^{(3)}\|_\infty \leq \frac{1}{8}M_0^3 M_3$, 
\begin{eqnarray*}\|k_3^{(3)}\|_\infty &\leq& 8M_0^3 M_3+\frac{27}{2}M_0^2 M_1M_2 + h(3M_0M_2 + \frac{49}{4}M_0M_1M_3+6M_1^2M_2)  \\ && +h^2(\frac{3}{2}M_0^3M_1M_2^2+6M_0^3M_1^2M_3) + h^3 M_0^3M_1^3M_3,
\end{eqnarray*}
where $M_i$ is an upper bound for $\|f^{(i)}\|_\infty$. For us, we have \[f(Y) = f(\rho,z,\alpha) = \left(\sin\alpha,\cos\alpha, \frac{z}{2}\sin\alpha + \left(\frac{1}{\rho}-\frac{\rho}{2}\right)\cos\alpha\right).\] On the region $\mathcal{R}$, we can then take $M_0 =2.671$, $M_1 = 9.34$, $M_2= 29.507$, $M_3=95.307$. 

Similarly, we may differentiate the equation $Y'(s)=f(Y(s))$ three times to estimate $\|Y^{(4)}\|_\infty$ in terms of the $M_i$, which yields $\|Y^{(4)}\|_\infty \leq 1532.075$. 

Combining these, we may take $C= 1229.881$. 

The global error for $s\in [0,\overline{s}]$ is then bounded by \[\frac{\exp(M_1\overline{s})}{M_1}Ch^3.\] We take $h=10^{-7}$. Substituting the values of $\overline{s},M_i,C$, we find that the global error is at most $(2.833 \times 10^{17}) h^3$. 
In particular since we chose $h=10^{-7}$ 
we obtain a global error of at most $\epsilon=0.0003$. 

That is, for $0\leq n \leq \overline{s}/h$, assuming that $Y_k \in \mathcal{R}$ for all $k<n$ and $Y(s) \in \mathcal{R}$ for all $s\in [0, nh]$, we have that \[ | Y(nh) - Y_n | \leq \epsilon = 0.0003,\] where $Y(s)$ is the exact solution. 

\begin{proof}[Proof of Proposition \ref{prop:numerical}]

As above, let $Y_n$ be the numerical solution using the RK3 method and $Y(s)$ the exact solution with initial conditions $Y(s) = (3.3,0,0)$. We will show by strong induction on $n$ that for all $n=0,\cdots, \overline{s}/h$, we have $Y_n\in \mathcal{R}$ and $Y(s)\in \mathcal{R}$ for all $0\leq s\leq \min((n+1)h,\overline{s})$.

Take $\delta=0.01$ and let 
\begin{equation}\label{eq:obsbounds}
\mathcal{R}_\delta = \{ (\rho,z,\alpha)| \rho \in [0.42 +\delta ,4-\delta] , |z|\leq 1-\delta ,  |\alpha| \leq \pi+1-\delta\}.
\end{equation}

Certainly the initial condition lies in $\mathcal{R}_\delta$, so since $M_0h \ll \delta$ we must have $Y(s)\in\mathcal{R}\delta$ for $s\in[0,h]$. Now at each step $n$, one may observe that the numerical solution lies in the region $\mathcal{R}_\delta$. (See Figures 1 and 2 
for a plot of the numerically obtained profile curve.) 
By the strong induction hypotheses and the error analysis in Section \ref{sec:error}, we have that $|Y(nh)-Y_n| \leq \epsilon \ll \delta$, hence $Y(nh) \in \mathcal{R}$. In fact, for any $s\in [(n-1)h,(n+1)h]$ we have $|Y(s)- Y_n| \leq \epsilon + M_0h  \ll \delta$. This implies that $Y(s) \in \mathcal{R}$ for all $s\in [0,(n+1)h]$, and completes the induction. 

In particular, the global error is at most $\epsilon$ for all $s\in[0,\overline{s}]$. For $s_1 = 3.745$, 
we then observe the particular value of the numerical solution:
\[ Y_\frac{s_1}{h} = (0.43648, 0.0076, -3.15049).\] 
Note that $3.15049 = \pi + 0.00889\cdots$. 
Since the global error is at most $\epsilon=0.0003$ we conclude that the exact solution $Y(s_1)$ satisfies $\alpha(s_1) < -\pi - 0.008 < -\pi$ but $z(s_1) > 0.007>0$ and hence $s_1 <s_2$. (See figures  3 to 8 
for plots comparing $\alpha$, $z$ and $s$.)
This completes the proof.
 \end{proof}

\begin{center}
\hspace{-1cm} $ 
\begin{matrix}\includegraphics[scale=0.65]{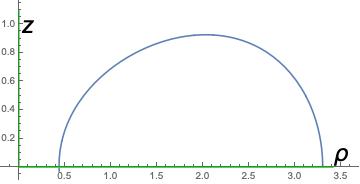}\\{Figure \ 1}
\end{matrix}
\begin{matrix}
\includegraphics[scale=0.6]{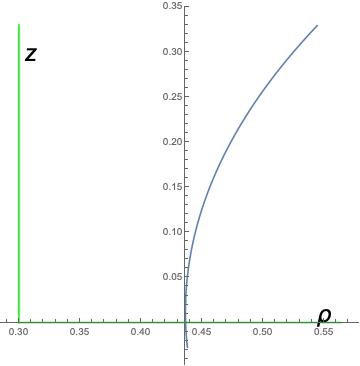}\\{Figure \ 2}
\end{matrix}
$
\end{center}

\begin{center}
 $ 
\begin{matrix}
\includegraphics[scale=0.6]{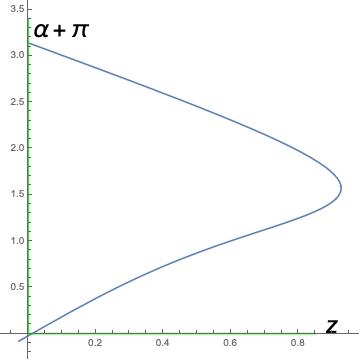}\\
 Figure \ 3 \end{matrix}
\begin{matrix}
\includegraphics[scale=0.6]{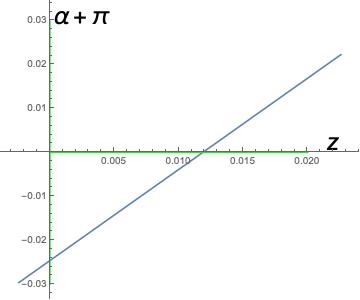}\\
Figure \ 4
\end{matrix}
$
\end{center}

\begin{center}
 $ 
\begin{matrix}
\includegraphics[scale=0.6]{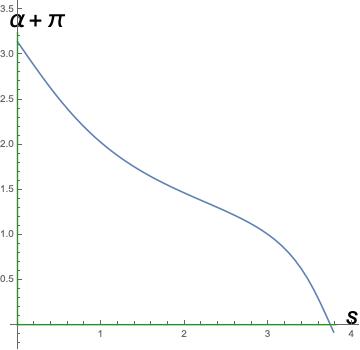}\\
 Figure \ 5 \end{matrix}
\begin{matrix}
\includegraphics[scale=0.6]{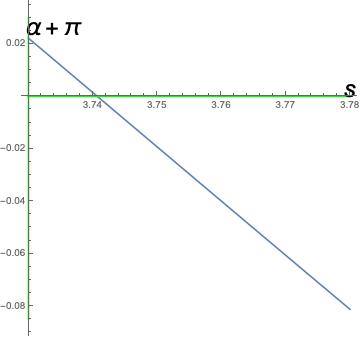}\\
 Figure \ 6 \end{matrix}
 $
\end{center}

\begin{center}
 $ 
\begin{matrix}
\includegraphics[scale=0.6]{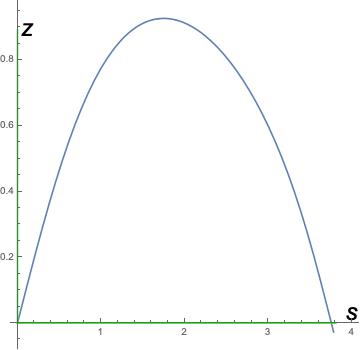}\\
 Figure \ 7 \end{matrix}
\begin{matrix}
\includegraphics[scale=0.6]{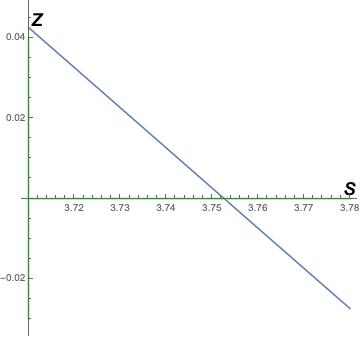}\\
 Figure \ 8 \end{matrix}
 $
\end{center}


\bibliographystyle{alpha}

{\small
\begin{thebibliography}{99}

\bibitem{An} Angenent, S.B.; Shrinking Doughnuts. {\it Nonlinear diffusion equations and the equilibrium states}, 3 (Gregynog 1989) 21-38 {\it Progr. Nonlinear differential equations appl.} 7, Birkhauser (1992).

\bibitem{BCP} Batista, M.; Cavalcante, M. P.; Pyo, J.; Some isoperimetric inequalities and eigenvalue estimates in weighted manifolds. {\it J. Math. Anal. Appl.} 419 (2014), no. 1, 617–626.

\bibitem{BeGo} Berger, Marcel; Gostiaux, Bernard; {\it Differential geometry: manifolds, curves, and surfaces}. GTM, 115. Springer-Verlag, New York, 1988.

\bibitem{BW}Bleecker, D., Weiner, J.: Extrinsic bounds on $\lambda_1$ of $\Delta$ on a compact manifold. {\it Comment. Math. Helv.} 51  (1976), 601-609.

\bibitem{br13} Brendle, Simon; Minimal surfaces in $S^3$: a survey on recent results. {\it Bull. Math. Sci.} 3 (2013), no.1, 133–171.

\bibitem{br15} Brendle, Simon;  
Embedded self-similar shrinkers of genus 0; {\em Ann. of Math.} 183 (2016), no. 2, 715–728

\bibitem {BK} Buser, P. and Karcher, H.; {\em
Gromov's almost flat manifolds},
 société mathematique de france, astérisque num. 81, Paris
 (1981).

\bibitem{CESG} Colbois, Bruno; El Soufi, Ahmad; Girouard, Alexandre; Isoperimetric control of the spectrum of a compact hypersurface. {\it J. Reine Angew. Math.} 683 (2013), 49–65.

\bibitem{Ch} Chavel, Isaac; {\it Eigenvalues in Riemannian geometry}. Pure and Applied Mathematics, 115. Academic Press, Inc., Orlando, FL, 1984. xiv+362 pp.

\bibitem{chmezh12b} Cheng, Xu;  Mejia, Tito and  Zhou, Detang; {Eigenvalue estimate and compactness for closed $f$-minimal surfaces},  {\it Pacific J. Math.} 271 (2014), no. 2, 347–367 
 
 \bibitem{CL} Coddington, Earl A.; Levinson, Norman; {\it Theory of ordinary differential equations}. McGraw-Hill Book Company, Inc., New York-Toronto-London, 1955. xii+429 pp.
 
\bibitem{comi} Colding, Tobias H.; Minicozzi, William P., II;  Generic mean curvature flow I: generic singularities. {\it Ann. of Math.} (2) 175 (2012), no. 2, 755–833.

\bibitem{dixi} Ding, Qi; Xin, Y. L. Volume growth, eigenvalue and compactness for self-shrinkers. {\it Asian J. Math.} 17 (2013), no. 3, 443–456.

\bibitem{DK} Drugan, G; Kleene, S.J; Inmersed self-shrinkers. {\em  Trans. Amer. Math. Soc.} 369 (2017), no. 10, 7213–7250.


\bibitem{DMWW} Du, Feng; Mao, Jing; Wang, Qiaoling; Wu, Chuanxi; { Eigenvalue inequalities for the buckling problem of the drifting Laplacian on Ricci solitons}; {\em Journal of Differential Equations}
 260, (2016) 5533-6450.
 
\bibitem{Ej} Ejiri, Norio;  A construction of nonflat, compact irreducible Riemannian manifolds which are isospectral but not isometric. {\it Math. Z.} 168 (1979), no. 3, 207–212.

\bibitem{GMO} Giménez, Fernando; Miquel, Vicente; Orengo, J. Javier;  Upper bounds of the first eigenvalue of closed hypersurfaces by the quotient area/volume. {\it Arch. Math}. 83 (2004), no. 3, 279--288.

\bibitem{HNW} Hairer, E; Norsett, S; Wanner, G; Solving Ordinary, Differential Equations I, Nonstiff problems, Vol.: 1. 2Ed. Springer-Verlag,(2000). 

\bibitem{he} Heintze, Ernst; Extrinsic upper bounds for $\lambda_1$. {\it Math. Ann.} 280 (1988), no. 3, 389–402.

\bibitem{huru} Huang, Qin; Ruan, Qihua;  Applications of some elliptic equations in Riemannian manifolds. {\it J. Math. Anal. Appl.} 409 (2014), no. 1, 189–196.

\bibitem{hui} Huisken, G.; Asymptotic behavior for singularities of the mean curvature flow, {\it J. Differential Geom. }31 (1990), 285–299. 


\bibitem{Ka} Karcher, H.;  Riemannian Center of Mass and Mollifier Smoothing {\em Comm. Pure and Applied Math.} {\bf XXX}, (1977) 509--541 

\bibitem{KM} Kleene, S. J; M\o ller, N. M; Self-shrinkers with a rotational symmetry. {\em Trans. Amer. Math. Soc.}, 366 (8): 3943-3963, 2014.  

\bibitem{lilo} Long, Li; On the spectrum of weighted Laplacian operator and its applications to uniqueness of Kähler Einstein metrics, {\it Math. Ann.} Published online 2014, DOI: 10.1007/s00208-014-1152-1

\bibitem{ma14} Ma, Li; Eigenvalue estimates and L1 energy on closed manifolds. {\it Acta Math. Sin.} (Engl. Ser.) 30 (2014), no. 10, 1729–1734.

\bibitem{madu} Ma, Li; Du, Sheng-Hua; Extension of Reilly formula with applications to eigenvalue estimates for drifting Laplacians. {\it C. R. Math. Acad. Sci. Paris} 348 (2010), no. 21-22, 1203–1206.

\bibitem{maliu} Ma, Li; Liu, Baiyu; Convexity of the first eigenfunction of the drifting Laplacian operator and its applications. {\it New York J. Math.} 14 (2008), 393–401.

\bibitem{mo11} Møller, Niels Martin ; Closed self-shrinking surfaces in $\mathbb{R}^3$ via the torus, Preprint, arXiv:1111.7318

\bibitem{piri} Pigola, Stefano; Rimoldi, Michele; Complete self-shrinkers confined into some regions of the space. {\it Ann. Global Anal. Geom.} 45 (2014), no. 1, 47–65

\bibitem{Re} Reilly, R.: On the first eigenvalue of the Laplacian for compact submanifold of
Euclidean space. {\it Comment. Math. Helv.} 52, 525-533 (1977)



\bibitem {sm05} Smoczyk, Knut;  Self-shrinkers of the mean curvature flow in arbitrary codimension. Int. Math. Res. Not. 2005, no. 48, 2983–3004.

\bibitem{sm12} Smoczyk, Knut; {\it Mean curvature flow in higher codimension: Introduction and survey}. in Global Differential Geometry, Springer Proceedings in Mathematics, vol. 17, part 2,(2012) 231-274. 


\end{thebibliography}

}

\bigskip

{\small 

\begin{tabular}{ l l }
&Departamento de Matem\'aticas para la Econom\'{\i}a y la Empresa \\
&Universidad de Valencia \\
&Avda Tarongers s/n\\
& 46022-Valencia (Spain)\\
      \       & email: carmen.domingo@uv.es\\
      
&  \\

&Departamento de Matemáticas\\
&Universidad de Valencia.\\
& Avda. Andr\'es Estell\'es, 1, 46100-Burjassot (Valencia) Spain \\
      \       & email: miquel@uv.es \\
      
      &  \\
      
&Department of Mathematics, \\
&Princeton University, \\
&Fine Hall, Washington Road, Princeton, NJ 08544, USA \\
      \       & email: jjzhu@math.princeton.edu 
     
\end{tabular}

}

\end{document}